# New Method for Constructing Complete Cap Sets

Iskandar Karapetyan[a, *], Karen Karapetyan[a, **]

[a]Department of Discrete Mathematics, Institute for Informatics and Automation Problems
of NAS RA, 0014, P. Sevak 1, Yerevan, Armenia

\* e-mail: isko@iiap.sci.am

\*\* e-mail: karen-karapetyan@iiap.sci.am

**Abstract**—A cap set in projective or affine geometry over a finite field is a set of points no three of which are collinear. In this paper, we propose a new construction for complete cap sets that yields a cap set of size 124928 in the affine geometry $AG(15,3)$. It should be noted that the constructed cap set in $AG(15,3)$ is more powerful and exceeds at least by 4096 points than those that can be obtained from the previously known ones using the product or doubling constructions.

## 1. INTRODUCTION

The cap set problem asks how large can a subset of affine geometry $AG(n,3)$ be that does not contain any line, that is, no three points $\boldsymbol{\alpha}, \boldsymbol{\beta}, \boldsymbol{\gamma}$ such that $\boldsymbol{\alpha} + \boldsymbol{\beta} + \boldsymbol{\gamma} = \mathbf{0}$ (*mod* 3). The problem was motivated by the design of statistical experiments (combinatorial design) [1], arithmetic progressions in the set of prime numbers [9, 10], coding theory [12], and a card game set [3], etc. More generally, the main problem in the theory of cap sets is to find the minimal and maximal sizes of complete cap sets in projective geometry $PG(n,q)$ or in affine geometry $AG(n,q)$. Note that the problem of determining the minimum size of a complete cap set in a given projective space is of particular interest in coding theory. If we write the points of the cap set as columns of a matrix, we obtain a matrix in which every three columns are linearly independent, hence the generator matrix of a linear orthogonal array of strength three. This matrix is a check matrix of a linear code with minimum distance greater than three. Many authors have noted that determining the exact value of the minimum and maximum cardinality of cap sets in the projective geometry $PG(n, q)$ or in the affine geometry $AG(n, q)$ seems to be a very hard, and each dimension n requires its own unique approach. There are some well-known constructions (product [4, 7] and doubling [17]), which allow to create large high-dimensional cap sets based on large low-dimensional cap sets. In this paper, we consider the problem of constructing complete cap sets in affine geometry $AG(n, 3)$ over the finite field $F_3 = \{0, 1, 2\}$. A cap set is called complete when it cannot be extended to a larger one. Let us denote the size of the largest cap set in $AG(n, q)$ and $PG(n, q)$ by $c_{n,q}$ and by $c'_{n,q}$, respectively. Presently, only the following exact values are known: $c_{n,2} = c'_{n,2} = 2^n$, $c_{2,q} = c'_{2,q} = q + 1$ if $q$ is odd, $c_{2,q} = c'_{2,q} = q + 2$ if $q$ is even, and $c_{3,q} = q^2$, $c'_{3,q} = q^2 + 1$ [1,20]. In addition to these general results, precise values of the maximum sizes of cap sets are known in the following cases: $c_{4,3} = c'_{4,3} = 20$ [8, 18], $c'_{5,3} = 56$ [11], $c_{5,3} = 45$ [6], $c_{4,4} = 40$, $c'_{4,4} = 41$ [5], $c_{6,3} = 112$ [19]. In the other cases, only lower and upper bounds on the sizes of cap sets in $AG(n, q)$ and $PG(n, q)$ are known [4,7,8]. In particular, it has been proven that $c_{7,3} \geq 236$ [2], using a computer search $c_{8,3} \geq 512$ [21]. In [15] we proved that $c_{11,3} \geq 5504$, and in [16] $c_{16,3} \geq 274432$, $c_{21,3} \geq 13991936$. Furthermore, using computer search, the following upper bounds on the maximum size of cap sets in dimensions seven to ten were proven in [22]: $c_{7,3} \leq 291$, $c_{8,3} \leq 771$, $c_{9,3} \leq 2070$, and $c_{10,3} \leq 5619$, respectively.

In this paper, we propose a new construction for complete cap sets, that yields a cap set of size 124928 in the affine geometry $AG(15,3)$. Note that the constructed cap set in $AG(15,3)$ is more powerful than those that can be obtained from the previously known ones using the product and doubling constructions.

## 2. NOTATIONS, DEFINITIONS AND KNOWN REZULTS

Note that three points $\alpha, \beta, \gamma$ in affine geometry $AG(n,3)$ are collinear if they are affinely dependent or, equivalently, $\alpha + \beta + \gamma = 0 \ (mod\ 3)$. Therefore, if $C_n$ is a cap set in $AG(n,3)$, then $\alpha + \beta + \gamma \neq 0 \ (mod\ 3)$ for every triple of distinct points $\alpha, \beta, \gamma \in C_n$. For the point $x \in AG(n,3)$, let's denote

$$x(0) = \{i \ | x_i = 0, i \in [1, n]\},$$

$$X(x) = \{\alpha | \alpha \in AG(n,3), \alpha(0) = x(0)\},$$

$$B_n = \{\alpha = (\alpha_1, \cdots, \alpha_n) \ | \alpha_i = 1, 2\}.$$

Obviously, $B_n$ is a cap set. Our new construction for cap sets is based on the concept of $P_n$-set, which was introduced in 2015 by K. Karapetyan [13]. A set of the points $A \subseteq AG(n,3)$ is called $P_n$-set, if it satisfies the following two conditions:

($i$) for any two distinct points $\alpha, \beta \in A$, there exists $i$ such that $\alpha_i = \beta_i = 0$, where $1 \leq i \leq n$,

($ii$) for any triple of distinct points $\alpha, \beta, \gamma \in A$, $\alpha + \beta + \gamma \neq 0 (mod\ 3)$.

Clearly, if the set $A \subseteq AG(n,3)$ is a $P_n$-set, then it is a cap set. To avoid new notations, in this paper we will simply denote the $P_n$-set as $P_n$ in $AG(n,3)$. We call $P_n$ to be odd, if $|\alpha(0)|$ is an odd number for every $\alpha \in AG(n, 3)$. $P_n$ is called complete when it cannot be extended to a larger one. The set $P_n$ is called $b$-saturated if $X(\alpha) \subseteq P_n$ for every point $\alpha \in P_n$, where $b = 1, 2$. Recall that concatenation or direct product of the sets is defined as follows. Let $A \subset AG(n, 3)$ and $B \subset AG(m, 3)$. Form a new set $AB \subset AG(n + m, 3)$ consisting of all points $\alpha = (\alpha_1, \cdots, \alpha_n, \alpha_{n+1}, \cdots, \alpha_{n+m})$, where $\alpha^1 = (\alpha_1, \cdots, \alpha_n) \in A$ and $\alpha^2 = (\alpha_{n+1}, \cdots, \alpha_{n+m}) \in B$. Concatenation for any number of sets is defined similarly. The next two theorems will introduce constructions to obtain a cap sets in higher dimensions by using known cap sets in lower dimensions. Theorem A states a simplified version of the general product construction theorem first stated by Mukhopadhyay in [17] and reformulated by Edel and Bierbrauer [4]. Theorem B, the doubling construction, is a special case of the general product construction. In this paper we will need the following theorems to construct new construction of complete cap sets.

Theorem A [4] (Product construction). Let $A \subseteq AG(n, 3)$ and $B \subseteq AG(m, 3)$ be cap sets. Then $AB \subset AG(n + m, 3)$ is a cap set.

Theorem B [17] (Doubling construction). Let $A \subseteq PG(n, 3)$ be a cap set. Then there is a cap set in $AG(n + 1, 3)$ of size $2|A|$.

Theorem C [15]. The set $A \subseteq AG(n, 3)$ is $b$-saturated and complete $P_n$-set if and only if it satisfies the following two conditions:

i) for any two distinct points $\alpha, \beta \in A$, there exists $i$ such that $\alpha_i = \beta_i = 0$, where $1 \leq i \leq n$,

iii) for any triple of distinct points $\alpha, \beta, \gamma \in A$, $\alpha(0) = \beta(0) = \gamma(0)$ or for two of them, say for $\alpha$ and $\beta$, there exists $i$ such that $\alpha_i = \beta_i = 0$ and $\gamma_i \neq 0$, where $1 \leq i \leq n$.

Theorem D [15]. If $P_n$ is $b$-saturated, complete and odd set, then $C_n = P_n \cup B'_n$ ($C_n = P_n \cup B''_n$) is a complete cap set.

For the given three sets $P_{n_1}$, $P_{n_2}$ and $P_{n_3}$, form the following set $P_{n_1}P_{n_2}B_{n_3} \cup P_{n_1}B_{n_2}P_{n_3} \cup B_{n_1}P_{n_2}P_{n_3}$. It is known that the formed set is a $P_n$-set [13], where $n = \sum_1^3 n_i$ and $n_1, n_2, n_3$ are any integers.

Theorem E («three» construction) [13]. The following recurrence relation $P_n = P_{n_1}P_{n_2}B_{n_3} \cup P_{n_1}B_{n_2}P_{n_3} \cup B_{n_1}P_{n_2}P_{n_3}$, with the initial sets $P_1 = \{(0)\}$, $P_2 = \{(0,1),(0,2)\}$ gives complete and $b$-saturated $P_n$ –set, where $n = \sum_{j=1}^3 n_j$.

For the given six sets $P_{n_1}, P_{n_2}, P_{n_3}, P_{n_4}, P_{n_5}$ and $P_{n_6}$, form the following ten sets:

$$A_1 = P_{n_1}P_{n_2}P_{n_3}B_{n_4}B_{n_5}B_{n_6},\ A_2 = P_{n_1}P_{n_2}B_{n_3}B_{n_4}B_{n_5}P_{n_6},$$

$$A_3 = P_{n_1}B_{n_2}P_{n_3}B_{n_4}P_{n_5}B_{n_6},\ A_4 = B_{n_1}P_{n_2}P_{n_3}P_{n_4}B_{n_5}B_{n_6},$$

$$A_5 = B_{n_1}B_{n_2}P_{n_3}P_{n_4}B_{n_5}P_{n_6},\ A_6 = B_{n_1}B_{n_2}P_{n_3}B_{n_4}P_{n_5}P_{n_6},$$

$$A_7 = B_{n_1}P_{n_2}B_{n_3}P_{n_4}P_{n_5}B_{n_6},\ A_8 = B_{n_1}P_{n_2}B_{n_3}B_{n_4}P_{n_5}P_{n_6},$$

$$A_9 = P_{n_1}B_{n_2}B_{n_3}P_{n_4}B_{n_5}P_{n_6},\ A_{10} = P_{n_1}B_{n_2}B_{n_3}P_{n_4}P_{n_5}B_{n_6}.$$

Theorem F («six» construction) [14]. The following recurrence relation $P_n = \bigcup_{i=1}^{10} A_i$, with the initial sets $P_1 = \{(0)\}$, $P_2 = \{(0,1),(0,2)\}$ gives complete and $b$-saturated $P_n$ -sets, where $n = \sum_{i=1}^6 n_i$ and $n_1, n_2, n_3, n_4, n_5, n_6$ are any integers.

# 3. MAIN RESULTS

Claim. If $x, y, z \in F_3$, then $x + y + z = 0 \pmod 3$ if and only if $x = y = z$ or they are pairwise distinct numbers.

Theorem. Suppose that the sets $P_n^1, P_n^2, P_n^3$ ($P_m^1, P_m^2, P_m^3$), and $P_k$ are complete, $b$-saturated and the following 3 conditions are satisfied:

1. If $x \in P_n^1, y \in P_n^2, z \in P_n^3$ ($x \in P_m^1, y \in P_m^2, z \in P_m^3$), then $x + y + z \neq 0 \pmod 3$.
2. If $x \in P_n^1, y, z \in P_n^3$ ($x \in P_m^1, y, z \in P_m^3$), then $x + y + z \neq 0 \pmod 3$.
3. If $x \in P_n^1 \cup P_n^2, y \in P_n^3$ ($x \in P_m^1 \cup P_m^2, y \in P_m^3$), then $x(0) \cap y(0) \neq \emptyset$.

Then the set $C = \bigcup_1^5 C_i$ is a complete cap set, where $k$ is any natural number and

$$C_1 = P_n^1 P_k B_m, C_2 = B_n P_k P_m^1, C_3 = P_n^2 B_k P_m^2, C_4 = P_n^3 B_k B_m, C_5 = B_n B_k P_m^3.$$

Proof. We want to prove that $C = \bigcup_1^5 C_i$ is a complete cap set. First, we prove that $C = \bigcup_1^5 C_i$ is a cap set in $AG(n + k + m, 3)$. Since every set $C_i$ is a product of a cap sets, therefore the Theorem A (product construction) implies that $C_i$ is a cap set, where $i \in \{1, 2, \dots, 5\}$. Thus, for each $i \in \{1, 2, \dots, 5\}$ there do not exist three distinct points $\alpha, \beta, \gamma \in C_i$ such that $\alpha + \beta + \gamma = 0 \pmod 3$. To prove that $C = \bigcup_1^5 C_i$ is a cap set, assume the opposite. Therefore, there exists a triple of distinct points $\alpha, \beta, \gamma \in C$ such that

$$\alpha + \beta + \gamma = 0 \pmod 3. \tag{*}$$

For convenience, we will further represent each point of $AG(n + k + m, 3)$ as a concatenation of three points. Thus, $\alpha = \alpha^1 \alpha^2 \alpha^3, \beta = \beta^1 \beta^2 \beta^3, \gamma = \gamma^1 \gamma^2 \gamma^3$, where $\alpha^1 = (\alpha_1, \cdots, \alpha_n), \alpha^2 = (\alpha_{n+1}, \cdots, \alpha_{n+k}), \alpha^3 = (\alpha_{n+k+1}, \cdots, \alpha_{n+k+m}), \beta^1 = (\beta_1, \cdots, \beta_n), \beta^2 = (\beta_{n+1}, \cdots, \beta_{n+k}), \beta^3 = $

$(\beta_{n+k+1}, \cdots, \beta_{n+k+m})$, $\boldsymbol{\gamma}^1 = (\gamma_1, \cdots, \gamma_n)$, $\boldsymbol{\gamma}^2 = (\gamma_{n+1}, \cdots, \gamma_{n+k})$, $\boldsymbol{\gamma}^3 = (\gamma_{n+k+1}, \cdots, \gamma_{n+k+m})$. We have 6 cases, depending on where $\boldsymbol{\alpha}, \boldsymbol{\beta}, \boldsymbol{\gamma}$ come from.

**Case $\boldsymbol{\alpha}, \boldsymbol{\beta}, \boldsymbol{\gamma} \in C_1 \cup C_2$:** Since, as mentioned above, each $C_i$ is a cap set ($i \in \{1, 2, \ldots, 5\}$), two points from the three points $\boldsymbol{\alpha}, \boldsymbol{\beta}, \boldsymbol{\gamma}$ contained in the same set, say $\boldsymbol{\alpha}, \boldsymbol{\beta} \in C_1$, and the third one is contained in $C_2$. From condition (i) of the $P_n^1$ set it follows that $\boldsymbol{\alpha}^1(0) \cap \boldsymbol{\beta}^1(0) \neq \emptyset$. Therefore, there exist $i$ such that $\alpha_i = \beta_i = 0$, where $i \in \{1, 2, \ldots, n\}$. Since $\boldsymbol{\gamma} \in C_2$, therefore $\boldsymbol{\gamma}^1 \in B_n$, hence $\gamma_i = 1$ or 2, and then the Claim implies that $\alpha_i + \beta_i + \gamma_i \neq 0$. Therefore, $\boldsymbol{\alpha} + \boldsymbol{\beta} + \boldsymbol{\gamma} \neq \mathbf{0} (mod\ 3)$, which contradicts equality (∗). If $\boldsymbol{\alpha}, \boldsymbol{\beta} \in C_2$ and $\boldsymbol{\gamma} \in C_1$, then the proof is similar to the one given above.

**Case $\boldsymbol{\alpha}, \boldsymbol{\beta} \in C_1 \cup C_2$, and $\boldsymbol{\gamma} \in \cup_3^5 C_i$:** According condition (i), for the set $P_k$, the following is true: $\boldsymbol{\alpha}^2(0) \cap \boldsymbol{\beta}^2(0) \neq \emptyset$. Therefore, there exist $i$ such that $\alpha_i = \beta_i = 0$, where $i \in \{n+1, \ldots, n+k\}$. Since $\boldsymbol{\gamma} \in \cup_3^5 C_i$, therefore $\boldsymbol{\gamma}^2 \in B_k$, hence $\gamma_i = 1$ or 2, and according to the Claim $\alpha_i + \beta_i + \gamma_i \neq 0$. Therefore, $\boldsymbol{\alpha} + \boldsymbol{\beta} + \boldsymbol{\gamma} \neq \mathbf{0} (mod\ 3)$, which again contradicts equality (∗). Therefore, no more than one of the points $\boldsymbol{\alpha}, \boldsymbol{\beta}, \boldsymbol{\gamma}$ can belong to $C_1 \cup C_2$, and the rest belong to $\cup_3^5 C_i$.

**Case $\boldsymbol{\alpha}, \boldsymbol{\beta} \in C_3$, and $\boldsymbol{\gamma} \in C \setminus C_3$:** Then, according to condition (i) of the set $P_n^2$ ($P_m^2$), the following holds: $\boldsymbol{\alpha}^1(0) \cap \boldsymbol{\beta}^1(0) \neq \emptyset$ ($\boldsymbol{\alpha}^3(0) \cap \boldsymbol{\beta}^3(0) \neq \emptyset$). Therefore, there exist $i$ such that $\alpha_i = \beta_i = 0$, where $i \in \{1, 2, \ldots, n\}$ ($i \in \{n+k+1, 2, \ldots, n+k+m\}$). Since $\boldsymbol{\gamma} \in C \setminus C_3$, therefore either $\boldsymbol{\gamma}^1 \in B_n$ or $\boldsymbol{\gamma}^3 \in B_m$, hence $\gamma_i = 1$ or 2, and it follows from the Claim that $\alpha_i + \beta_i + \gamma_i \neq 0$ for some $i$, where $i \in \{1, 2, \ldots, n\}$ or $i \in \{n+k+1, 2, \ldots, n+k+m\}$. Therefore, $\boldsymbol{\alpha} + \boldsymbol{\beta} + \boldsymbol{\gamma} \neq \mathbf{0} (mod\ 3)$, which contradicts equality (∗). Therefore, no more than one of the points $\boldsymbol{\alpha}, \boldsymbol{\beta}, \boldsymbol{\gamma}$ can belong to $C_3$, and the rest belong to $C \setminus C_3$.

**Case $\alpha \in C_1 \cup C_2, \beta \in C_3$ and $\gamma \in C_4 \cup C_5$:** If $\alpha \in C_1$ and $\gamma \in C_4$ ($\alpha \in C_2$ and $\gamma \in C_5$), then by condition (1) of the Theorem $\alpha^1 + \beta^1 + \gamma^1 \neq 0 (mod\ 3)(\alpha^3 + \beta^3 + \gamma^3 \neq 0 (mod\ 3))$. Therefore, $\alpha + \beta + \gamma \neq 0 (mod\ 3)$. If, $\alpha \in C_2 (\alpha \in C_1)$ and $\gamma \in C_4 (\gamma \in C_5)$, then by condition (3) of the Theorem $\gamma^1(0) \cap \beta^1(0) \neq \emptyset$ ($\beta^3(0) \cap \gamma^3(0) \neq \emptyset$). Therefore, there exist $i$ such that $\gamma_i = \beta_i = 0$, where $i \in \{1, \dots, n\} (i \in \{n+k+1, \dots, n+k+m\})$. Since $\alpha^1 \in B_n (\alpha^3 \in B_m)$, hence $\alpha_i = 1$ or 2, and it follows from the Claim that $\alpha_i + \beta_i + \gamma_i \neq 0$, where $i \in \{1, \dots, n\} (i \in \{n+k+1, \dots, n+k+m\})$. Therefore, $\alpha + \beta + \gamma \neq 0 (mod\ 3)$, which again contradicts equality (∗).

**Case $\alpha \in C_1 \cup C_2, \beta, \gamma \in C_4 \cup C_5$:** If $\alpha \in C_1$ ($\alpha \in C_2$) and $\beta, \gamma \in C_4 (\beta, \gamma \in C_5)$, then according condition (2) of the Theorem $\alpha^1 + \beta^1 + \gamma^1 \neq 0(mod\ 3)(\alpha^3 + \beta^3 + \gamma^3 \neq 0(mod\ 3))$. Therefore, $\alpha + \beta + \gamma \neq 0(mod\ 3)$, which again contradicts equality (∗). If $\alpha \in C_1$ ($\alpha \in C_2$) and $\beta, \gamma \in C_5 (\beta, \gamma \in C_4)$, then according condition (i) of the $P_m^3 (P_n^3)$ set $\beta^3(0) \cap \gamma^3(0) \neq \emptyset$ ($\beta^1(0) \cap \gamma^1(0) \neq \emptyset$). Therefore, there exist $i$ such that $\gamma_i = \beta_i = 0$, where $i \in \{n+k+1, \dots, n+k+m\} (i \in \{1, \dots, n\})$. Since $\alpha^3 \in B_m (\alpha^1 \in B_n)$, hence $\alpha_i = 1$ or 2, and the Claim implies that $\alpha_i + \beta_i + \gamma_i \neq 0$, where $i \in \{n+k+1, \dots, n+k+m\} (i \in \{1, \dots, n\})$. Therefore, $\alpha^3 + \beta^3 + \gamma^3 \neq 0(mod\ 3)(\alpha^1 + \beta^1 + \gamma^1 \neq 0\ (mod\ 3))$, contradicting equality (∗). If $\alpha \in C_1$ ($\alpha \in C_2$) and $\beta \in C_4, \gamma \in C_5$, then according condition (3) of the Theorem $\alpha^1(0) \cap \beta^1(0) \neq \emptyset$ ($\alpha^3(0) \cap \gamma^3(0) \neq \emptyset$). Since $\gamma^1 \in B_n (\beta^3 \in B_m)$, hence $\alpha_i = 1$ or 2, and the Claim follows that $\alpha_i + \beta_i + \gamma_i \neq 0$, where $i \in \{1, \dots, n\} (i \in \{n+k+1, \dots, n+k+m\})$. Therefore, $\alpha^1 + \beta^1 + \gamma^1 \neq 0(mod\ 3)(\alpha^3 + \beta^3 + \gamma^3 \neq 0(mod\ 3))$. Therefore, $\alpha + \beta + \gamma \neq 0(mod\ 3)$, which again contradicts equality (∗).

**Case $\alpha, \beta, \gamma \in A_4 \cup A_5$:** Since $C_4$ and $C_5$, as mentioned above, are cap sets, therefore two points from the three points $\alpha, \beta, \gamma$ contained in the same set, say $\alpha, \beta \in C_4$, and the third one is contained in $C_5$. Then according condition (i) of the $P_n^3$, $\alpha^1(0) \cap \beta^1(0) \neq \emptyset$. Since $\gamma^1 \in B_n$, hence $\gamma_i = 1$ or 2, and the Claim implies that $\alpha_i + \beta_i + \gamma_i \neq 0$, where $i \in \{1, ..., n\}$. Therefore, $\alpha + \beta + \gamma \neq 0 \pmod{3}$, contradicting equality (∗). If $\alpha, \beta \in C_5$ and $\gamma \in C_4$, then the proof is similar to the one above.

Now we will prove the completeness of the set $C = \cup_1^5 C_i$ by contradiction. Suppose that $C = \cup_1^5 C_i$ is not a complete cap set. Therefor there is a point $x \notin C$ such that $\{x\} \cup C$ is a cap set. Let`s represent the point $x$, as shown above, by concatenation of three points. Thus, $x = x^1 x^2 x^3$, where $x^1 = (x_1, \cdots, x_n)$, $x^2 = (x_{n+1}, \cdots, x_{n+k})$, $x^3 = (x_{n+k+1}, \cdots, x_{n+k+m})$. Then $\alpha^1(0) \cap x^1(0) \neq \emptyset$ for each $\alpha^1 \in P_n^3$ or there exits $\alpha^1 \in P_n^3$ such that $\alpha^1(0) \cap x^1(0) = \emptyset$. In the first case, the completeness of the set $P_n^3$ implies that $x^1 \in P_n^3$, or there are two points $y^1, z^1 \in P_n^3$, such that $x^1 + y^1 + z^1 = 0 \pmod{3}$. Then we can choose $y^2 = (y_{n+1}, \cdots, y_{n+k})$, $z^2 = (z_{n+1}, \cdots, z_{n+k}) \in B_k$, $y^3 = (y_{n+k+1}, \cdots, y_{n+k+m})$, $z^3 = (z_{n+k+1}, \cdots, z_{n+k+m}) \in B_m$ as follows: If $x_i = 0$, then we take $y_i$ and $z_i$ in such a way that each of them is the additive complement of the other, otherwise, we take $y_i = z_i = x_i$, where $i \in \{n+1, ..., n+k+m\}$. Then $y = x^1 y^2 y^3$, $z = x^1 z^2 z^3 \in C_4$, or $y = y^1 y^2 y^3$, $z = z^1 z^2 z^3 \in C_4$ and in both subcases $x + y + z = 0 \pmod{3}$, which contradicts the assumption that $\{x\} \cup C$ is a cap set. Therefore, there exits $\alpha^1 \in P_n^3$ such that $\alpha^1(0) \cap x^1(0) = \emptyset$. By the same reason, for the point $x^3$ there exits point $\gamma^3 \in P_m^3$ such that $\gamma^3(0) \cap x^3(0) = \emptyset$. The b-saturatness of the $P_n^3$ implies that $X(\alpha^1) \subseteq P_n^3$. We will choose two points $y^1$ and $z^1$ in the following way: If $x_i = 0$, then we take $z_i$ as the additive complement of the $\alpha_i$. If $\alpha_i = 0$, then we take $z_i$ as the additive complement of the $x_i$, otherwise we take $y_i = z_i = x_i$, where $i \in \{1, ..., n\}$. Obviously $y^1 \in P_n^3$, $z^1 \in B_n$ and $x^1 + y^1 + z^1 = 0 \pmod{3}$. Similarly, one can choose $y^3 \in P_m^3$ and $z^3 \in B_m$ such that $x^3 + y^3 + z^3 = 0 \pmod{3}$. Now we will choose $y^2$ and $z^2 \in B_k$ as follows: If $x_i = 0$, then we take $y_i$ and $z_i$ in such a way that each of them is the additive complement of the other, otherwise, we take $y_i = z_i = x_i$, where $i \in \{n+1, \ldots, n+k\}$. Clearly, the choice of the points $y^2$ and $z^2$ implies that $x^2 + y^2 + z^2 = 0 \pmod{3}$. Obviously, the last three equalities imply that $x + y + z = 0 \pmod{3}$, which contradicts our supposition that $\{x\} \cup C$ is a cap set. Theorem is proved.

Corollary. $c_{15,3} \geq 124928$.

Proof. Let`s take $k = 3$ and $n = m = 6$. Theorem E for $n_1 = n_2 = n_3 = 1$ implies that $P_3 = \{(0,0,1), (0,0,2), (0,1,0), (0,2,0), (1,0,0), (2,0,0)\}$. Theorem F for $n_i = 1$ implies that $P_6^1 = \bigcup_{i=1}^{10} A_i$, where

$$A_1 = P_1 P_1 P_1 B_1 B_1 B_1, \; A_2 = P_1 P_1 B_1 B_1 B_1 P_1,$$

$$A_3 = P_1 B_1 P_1 B_1 P_1 B_1, \; A_4 = B_1 P_1 P_1 P_1 B_1 B_1,$$

$$A_5 = B_1 B_1 P_1 P_1 B_1 P_1, \; A_6 = B_1 B_1 P_1 B_1 P_1 P_1,$$

$$A_7 = B_1 P_1 B_1 P_1 P_1 B_1, \; A_8 = B_1 P_1 B_1 B_1 P_1 P_1,$$

$$A_9 = P_1 B_1 B_1 P_1 B_1 P_1, \; A_{10} = P_1 B_1 B_1 P_1 P_1 B_1.$$

$P_1 = \{(0)\}$, $B_1 = \{(1), (2)\}$ and $i \in \{1, \ldots, 6\}$. Let`s denote the mirror inversion of the set $P_6^1$ by $P_6^2$. It is easy to check that $P_6^2 = \bigcup_{1}^{10} A_i'$, where

$$A_1' = B_1 B_1 B_1 P_1 P_1 P_1, \; A_2' = P_1 B_1 B_1 B_1 P_1 P_1,$$

$$A_3' = B_1 P_1 B_1 P_1 B_1 P_1, \; A_4' = B_1 B_1 P_1 P_1 P_1 B_1,$$

$$A_5' = P_1 B_1 P_1 P_1 B_1 B_1, \; A_6' = P_1 P_1 B_1 P_1 B_1 B_1,$$

$$A_7' = B_1 P_1 P_1 B_1 P_1 B_1, \; A_8' = P_1 P_1 B_1 B_1 P_1 B_1,$$

$$A'_9 = P_1B_1P_1B_1B_1P_1, \quad A'_{10} = B_1P_1P_1B_1B_1P_1.$$

As in [2], we take $P_6^3 = \{e_{1,1} = (1, 0, ..., 0), e_{1,2} = (2, 0, ..., 0), e_{2,1} = (0, 1, ..., 0), e_{2,2} = (0, 2, ..., 0), ..., e_{6,1} = (0, 0, ..., 1), e_{6,2} = (0, 0, ..., 2)\} \subset AG(6, 3)$. It is obvious that $|P_6^3| = 12$. Clearly, the set $P_6^3$ is $P_6$-set. Since for every pair of points $\boldsymbol{\alpha} \in P_6^1$ and $\boldsymbol{\beta} \in P_6^2$, $\boldsymbol{\alpha}(0) \nsubseteq \boldsymbol{\beta}(0)$ and $\boldsymbol{\beta}(0) \nsubseteq \boldsymbol{\alpha}(0)$. Therefore, there are $i, j, i \neq j$ such that $\alpha_i = 0, \beta_i \neq 0$ and $\alpha_j \neq 0, \beta_j = 0, i, j \in \{1, 2, ..., 6\}$. Since each point $\boldsymbol{\gamma} \in P_6^3$ has only one non-zero coordinate, therefore it follows from the Claim that $\alpha_i + \beta_i + \gamma_i \neq 0 (mod\ 3)$ or $\alpha_j + \beta_j + \gamma_j \neq 0 (mod\ 3)$, hence $\boldsymbol{\alpha} + \boldsymbol{\beta} + \boldsymbol{\gamma} \neq \boldsymbol{0}(mod\ 3)$. Thus, the condition 1 of the Theorem is satisfied. Since each point $\alpha \in P_6^1 (\alpha \in P_6^2)$ has three non-zero coordinates, but every point $\boldsymbol{\gamma} \in P_6^3$ has only one non-zero coordinate, therefore, for each $\boldsymbol{\alpha} \in P_6^1 (\boldsymbol{\alpha} \in P_6^2)$ and each pair $\boldsymbol{\beta}, \boldsymbol{\gamma} \in P_6^3$, $\boldsymbol{\alpha} + \boldsymbol{\beta} + \boldsymbol{\gamma} \neq \boldsymbol{0}(mod\ 3)$. Thus, the condition 2 of the Theorem is also satisfied. The condition 3 is obviously satisfied, since each point $\boldsymbol{\alpha} \in P_6^1 \cup P_6^2$ has three zero coordinates, and each point $\boldsymbol{\gamma} \in P_6^3$ has five zero coordinates, therefore, $\boldsymbol{\alpha}(0) \cap \boldsymbol{\beta}(0) \neq \emptyset$. Now, replacing the sets $P_n^1, P_n^2, P_n^3 (P_m^1, P_m^2, P_m^3)$ in each $C_i$ with $P_6^1, P_6^2, P_6^3$ respectively, then we get a set $C = \bigcup_1^5 C_i$ in $AG(15, 3)$, where $i \in \{1, 2, ..., 5\}$. Clearly, $C_i \cap C_j = \emptyset$ for every pair $i, j, i \neq j$, where, $i, j \in \{1, 2, ..., 5\}$. Thus, $|C| = \sum_1^5 |C_i|$. It is easy to check that $|P_6^1| = |P_6^2| = 80$, and $|B_6| = 64$. Then it follows from the Theorem that $|C| = 2 * 80 * 6 * 64 + 80 * 8 * 80 + 2 * 12 * 8 * 64 = 124928$. Corollary is proved.

## References

1. R. C. Bose, "Mathematical theory of the symmetrical factorial design", *Sankhya*, vol. 8, pp. 107-166, 1947, MR0026781.

2. A. R. Calderbank and P. C. Fishburn, "Maximal three-independent subsets of $\{0, 1, 2\}^n$", *Designs, Codes and Cryptography,* vol. 4, pp. 203-211, 1994, DOI 10.1007/BF01388452, MR1277940.

3. B. L. Davis and D. Maclagan, "The card game set", *The Mathematical Intelligencer,* vol. 25(3), pp. 33-40, 2003, DOI 10.1007/BF0298484, MR2005098.

4. Y. Edel, and J. Bierbrauer, "Recursive constructions for large caps", *Bulletin.of the Belgian Math. Society* - Simon Stevin. 6, pp. 249-258, 1999, DOI 10.36045/bbms/1103141034.

5. Y. Edel and J. Bierbrauer, "41 is the largest size of a cap in $PG(4,4)$", *Designs, Codes and Cryptography*, vol. 16, pp. 151-160, 1999, DOI 10.1023/A:1008389013117.

6. Y. Edel, S. Ferret, I. Landjev and L. Storme, "The classification of the largest caps in $AG(5,3)$", *Journal of Combinatorial Theory*, ser. A, vol. 99, pp. 95-110, 2002, DOI 10.1006/jcta.2002.3261, MR1911459.

7. Y. Edel, "Extentions of generalized product caps", *Designs, Codes and Cryptography*, vol. 31, pp. 5-14, 2004, DOI 10.1023/A:1027365901231, MR2031694.

8. Y. Edel and J. Bierbrauer, "Large caps in small spaces", *Designs, Codes and Cryptography*, vol. 23(2), pp. 197-212, 2001, DOI 10.1023/A:1011216716700.

9. J. S. Ellenberg and D. Gijswijt, "On large subsets of $F_q^n$ with no three-terms arithmetc progression", *Ann. Of Math*. (2) 185 (2017), no. 1, pp. 339-343, DOI 10.4007/annals.2017.185.1.8, MR3583358.

10. JOSHA A. GROCHOW, "New applications of the polynomial method: The cap set conjecture and beyond", *BULLETIN (New Series) OF THE AMERICAN MATHEMATICAL SOCIETY*, Vol. 56, Num. 1, January, pp. 29-64, 2019, DOI 10.1090/bull/1648.

11. R. Hill, "On the largest size of cap in $S_{5,3}$", *Atti Accad. Naz.* Lincei *Rendiconti*, vol. 54, pp. 378-384, 1973.

12. J. W. P. Hirschfeld and L. Storme, "The packing problem in statistics, coding theory and finite projective spaces", update 2001, in *Finite geometry*, vol. 3 of Dev. Math., Kluwer Acad. Publ., Dordrecht, pp. 201-246, 2001, DOI 10.1007/978-1-4613-0283-4_13.

13. K. Karapetyan, "Large Caps in Affine Space $AG(n,3)$", *Proceedings of International Conference Computer Science and Information Technologies*, Yerevan, Armenia, pp. 82-83, 2015.

14. K. Karapetyan, "On the complete caps in Galois affine space $AG(n,3)$", *Proceedings of International Conference Computer Science and Information Technologies*, Yerevan, Armenia, p. 205, 2017.

15. Karapetyan and K. Karapetyan, "Complete caps in affine geometry AG(n, 3)," Patern Recognit Image Anal, 34, 74–91 (2024). https://doi.org/10.1134/s1054661824010097

16. Iskandar Karapetyan and Karen Karapetyan, "On Some Large Cap Sets", Pattern Recognition and Image Analysis, 2025, Vol. 35, No. 4, ISSN 1054-6618, pp. 884–888, DOI 10.1134/S1054661825700683.

17. A. C. Mukhopadhyay, "Lower bounds on $m_t(r,s)$", *Journal of Combinatorial Theory,* Seires A, vol. 25(1), pp. 1-13, 1971.

18. G. Pellegrino, "Sul Massimo ordine delle calotte in $S_{4,3}$", *Matematiche* (Catania), vol. 25, pp. 149-157, 1971, MR0363652.

19. A. Potechin, "Maximal caps in $AG(6,3)$", *Designs, Codes and Cryptography*, vol. 46, pp. 243-259, 2008, DOI 10.1007/s10623-007-9132-z.

20. B. Qvist, "Some remarks concerning curves of the second degree in a finite plane", *Ann Acad. Sci. Fenn*, Ser. A, vol. 134, p. 27. 1952.

21. B. Romera-Paredes, M. Barekatain, A. Novikov, M. Balog, M. P. Kumar, E. Dupont, F. J. R. Ruiz, J. S. Ellenberg, P. Wang, O. Fawzi, P. Kohli, and A. Fawzi, "Mathematical discoveries from program search with large language models," Nature 625, 468–475 (2024). https://doi.org/10.1038/s41586-023-06924-6

22. N. D. Versluis, "On The Cap Set Problem, Upper bounds on maximal cardinalities of caps in dimensions seven to ten", Delft University of Technology, pp. 1-52, July 2017.